\documentclass[notitlepage,leqno,10pt]{article}

\usepackage[latin1]{inputenc}
\usepackage{amsfonts}
\usepackage{amsmath}
\usepackage{amssymb}
\usepackage{amsrefs}

\usepackage{color}
\definecolor{green2}{rgb}{0,0.6,0}

\newcommand{\R}{\mathbb{R}}
\newcommand{\N}{\mathbb{N}}

\newcommand{\qed}{\penalty 500\hfill$\square$\par\medskip}

\newcommand{\dist}{{\rm dist}}

\newcommand{\eps}{\varepsilon}

\def\XXint#1#2#3{{\setbox0=\hbox{$#1{#2#3}{\int}$}
     \vcenter{\hbox{$#2#3$}}\kern-.5\wd0}}

\def\sub{\underline}
\def\super{\overline}

\newtheorem{lem}         {Lemma}[section]

\newtheorem{pro}    [lem]{Proposition}
\newtheorem{thm}    [lem]{Theorem}
\newtheorem{rem}    [lem]{Remark}

\title{Uniqueness of large solutions}
\date{\today}
\author{O. Costin$^{1}$, L. Dupaigne$^{2}$ and O. Goubet$^{2}$}
\begin{document}
\maketitle

{\begin{center}
$^{1}${\small
Department of Mathematics, The Ohio State University,\\
100 Math Tower, 231 West 18th Avenue, Columbus, OH 43210-1174, USA\\ {\tt costin@math.ohio-state.edu}\\
\smallskip
$^2$LAMFA, UMR CNRS 7352, Universit\'e Picardie Jules Verne \\33, rue St Leu, 80039 Amiens, France \\  \tt louis.dupaigne@math.cnrs.fr, olivier.goubet@u-picardie.fr\\
}
\end{center}}

\abstract{Given a nondecreasing nonlinearity $f$, we prove uniqueness of large solutions to the equation \eqref{main} below, in the following two cases: the domain is the ball or the domain has nonnegative mean curvature and the nonlinearity is asymptotically convex.}

\section{Introduction}
In this paper, we are interested in the so-called large solutions of a certain class of partial differential equations. Let us recall what they are:
given $\Omega$ be a bounded domain of  $\R^N$, $N\ge1$ and $f\in C^1(\R)$, a large solution is a function $u\in C^2(\Omega)$ satisfying
\begin{equation}\label{main} 
 \left\{
 \begin{aligned} 
\Delta u &=f(u) &\quad\text{in $\Omega$,}\\
u &= +\infty &\quad\text{on $\partial \Omega$,}
\end{aligned}
\right. 
 \end{equation} 
where the boundary condition is understood in the sense that
$$
\lim_{x\to x_{0}, x\in \Omega}u(x) = +\infty\qquad\text{for all $x_{0}\in\partial \Omega$}
$$
and where $f$ is assumed to be positive at infinity, in the sense that
\begin{equation} \label{positivity} 
\exists\; a\in\R\quad\text{s.t.}\quad f(a)>0\quad\text{and} \quad f(t)\ge0\quad\text{for $t>a$.}
\end{equation} 
When the boundary of $\Omega$ is smooth enough, existence of a solution of \eqref{main} is equivalent to the so-called Keller-Osserman condition~: 
\begin{equation}\label{KO}
\int^{+\infty}\frac{dt}{\sqrt{F(t)}}<+\infty,\qquad\text{where $F(t)=\int_{a}^{t}f(s)\;ds$.}
\end{equation} 
For a proof of this fact, see the seminal works of J.B. Keller \cite{keller} and R. Osserman \cite{osserman} for the case of monotone $f$, as well as \cite{ddgr} for the general case. From here on, we always assume that \eqref{KO} holds. 

Uniqueness of solutions of \eqref{main} turns out to be delicate. As one might expect, it fails in the presence of oscillations. For example, if $f(u)=u^2\sin^2(u)$, the equation has infinitely many solutions (see \cite{ddgr}). It is also known (see e.g. the remark on p. 325 in \cite{veron3}) that uniqueness fails for a nonlinearity of the form $f(u)=u^p$, $p>1$, if the domain is not smooth enough: 
\begin{pro}\label{pro1}Assume that $\Omega=B\setminus\{0\}$ is the punctured unit ball of $\R^N$, $N\ge2$. Let $p\in(1,\frac N{N-2})$ if $N\ge 3$ (respectively $p\in(1,+\infty)$ if $N=2$) and $f(u)=u^p$. Then, there exists infinitely many solutions of \eqref{main}. 
\end{pro}
However, one could hope that uniqueness holds under the simple assumptions that $f$ is a nondecreasing function and that $\Omega$ has smooth boundary. As of today, this question remains open.
In \cite{cd}, we proved uniqueness in the case where $\Omega$ is a ball. 
\begin{thm}[\cite{cd}]\label{thm1} Assume that $\Omega$ is the unit ball in $\R^N$, $N\ge1$. Assume that $f$ is a nondecreasing function such that \eqref{positivity} and \eqref{KO} hold. Then, there exists a unique solution of \eqref{main}. 
\end{thm}
In this paper, we give a shorter proof of this fact. Under extra convexity assumptions, we obtain the following answer for a more general class of domains.   
\begin{thm}\label{thm2} Assume that $\partial\Omega$ is of class $C^3$ and that its mean curvature is nonnegative. Assume that $f$ is a nondecreasing function such that \eqref{positivity} and \eqref{KO} hold. Assume in addition that there exists $M\in\R$ such that $\sqrt{F}$ is convex in $(M,+\infty)$. Then, there exists a unique solution of \eqref{main}. 
\end{thm}
\begin{rem}If $f$ is asymptotically convex, then so is $\sqrt F$.
\end{rem}

Let us turn to the proofs.
\section{Proof of Theorem \ref{thm1}}

\noindent {\bf Step 1.} Reduction to the radial case.

Assume $\Omega$ is the ball.
It is well-known (see e.g. Lemma 2.4 in \cite{cd}) that the equation has a minimal and a maximal solution, each of which is radial. That is, there exist two large radial solutions $U_{1}, U_{2}$ such that any large solution $u$ satisfies $U_{1}\le u\le U_{2}$. In particular, it suffices to prove that $U_{1}\ge U_{2}$.

\noindent {\bf Step 2.} Let $u$ be a large radial solution. There exists $r_{0}\in(0,1)$ such that in $(r_{0},1)$, $u$ is strictly increasing and
\begin{equation} \label{keller}
\frac1{2N}F(u)\le \left(\frac{du}{dr}\right)^2\le 4 F(u)
\end{equation}  
This is essentially Keller's classical argument (see \cite{keller}): let $u$ be a large radial solution. Using \eqref{positivity}, it follows that for $r$ close to $1$,
\begin{equation} \label{rad} 
r^{1-N}\frac{d}{dr}\left(r^{N-1}\frac{du}{dr}\right)=\Delta u=f(u)\ge0.
\end{equation} 
Since $u$ is unbounded, there exists $r_{1}$ close to $1$ such that $du/dr(r_{1})>0$. By \eqref{rad},  $du/dr>0$ in $[r_{1},1)$. Integrating \eqref{rad}, we also have for $r\in (r_{1},1)$,
\begin{align*}
r^{N-1}\frac{du}{dr} &= r_{1}^{N-1}\frac{du}{dr}(r_{0}) + \int_{r_{1}}^r s^{N-1}f(u(s))\;ds \\
&\le r_{1}^{N-1}\frac{du}{dr}(r_{1})+  f(u(r))\frac{r^N}{N}. 
\end{align*} 
Since $f$ is nondecreasing and satisfies the Keller-Osserman condition \eqref{KO}, $\lim_{+\infty}f=+\infty$.  Using this in the above, given $\epsilon>0$, we find $r_{2}\in [r_{1},1)$ such that for  $r\in (r_{2},1)$,
$$
\frac{1}r\frac{du}{dr}\le \left(\frac1N+\epsilon\right) f(u)
$$
Taking $\epsilon=\frac1{2(N-1)}$ and recalling that 
$$
\frac{d^2u}{dr^2}+\frac{N-1}{r}\frac{du}{dr} = f(u),
$$
we deduce that
$$
\frac1{2N}f(u)\le\frac{d^2u}{dr^2}\le f(u)\qquad\text{in $[r_{2},1)$.}
$$
Multiplying by $2du/dr$, integrating and letting $c=du/dr(r_{2})^2 - F(u(r_{2}))$, we obtain
$$
\frac1N F(u) +c \le \left(\frac{du}{dr}\right)^2 \le 2 F(u) + c \qquad\text{for $r\in[r_{2},1)$}
$$
and so we find $r_{0}\in[r_{2},1)$ such that \eqref{keller} holds in $[r_{0},1)$. 

\noindent{\bf Step 3.} Change of independent variable.

Thanks to Step 2, for $r$ close to 1, given $i\in\{1,2\}$, we may perform the change of variable $u=U_{i}(r)$. Let $r=r_{i}(u)$ denote the inverse mapping of $U_{i}$ and $V_{i}=\frac{dU_{i}}{dr}\circ r_{i}$. By the chain rule,
\begin{equation} \label{chain rule}  V_{i}\frac{dV_{i}}{du}+\frac{N-1}{r_{i}}V_{i}=f, \end{equation}
while $dr_{i}/du=1/V_{i}$, so that
\begin{equation} \label{r equation}  1-r_{i}=\int_u^{+\infty} \frac 1 {V_{i}} du'.\end{equation}

\noindent {\bf Step 4.} There exists $u_{0}>0$  such that $r_1\ge r_2$ and $V_1\ge V_2$ in $[u_{0},+\infty)$.

Since $r_{i}$ is the inverse mapping of $U_{i}$ and $U_{1}\le U_{2}$, we have $r_{1}\ge r_{2}$. By \eqref{chain rule}, the function $z=V_{2}-V_{1}$ satisfies
$$
\frac{dz}{du}+(N-1)\left\{\frac1{r_{2}}-\frac1{r_{1}}\right\} = \left(\frac1{V_{2}}-\frac1{V_{1}}\right)f = -\frac{f}{V_{1}V_{2}}z.
$$
Since $r_{1}\ge r_{2}$, we deduce that $w$ satisfies the differential inequality
\begin{equation} \label{gronwall}
\frac{dz}{du}+a z\le 0,
\end{equation}
where $a= \frac{f}{V_{1}V_{2}}\ge0$ for large $u$.
By \eqref{r equation}, we also have
$$
\int_{u}^{+\infty}\frac{1}{V_{2}}du' \ge \int_{u}^{+\infty}\frac1{V_{1}}du'.
$$
So, there must exist $u_{0}$ such that $1/V_{2}(u_{0})\ge 1/V_{1}(u_{0})$ i.e. $w(u_{0})\le 0$. Using this together with \eqref{gronwall}, we deduce that $z\le0$ in $[u_{0},+\infty)$, as desired.

\noindent {\bf Step 5.} { The function $w=r_1^{2N-2}V_1^2-r_2^{2N-2}V_2^2$ is bounded. 

To see this, observe first that

\begin{equation} \label{step5}  \frac{dw}{du}=2(r_1^{2N-2}-r_2^{2N-2})f.\end{equation} 

\noindent Hence, $w$ is a nonnegative nondecreasing function and
$$\frac{dw}{du}\le 4(N-1)(r_{1}-r_{2})f=4(N-1)\left(\int_{u}^{+\infty}\left(\frac{1}{V_2}-\frac{1}{V_1}\right)du'\right)f\\
$$
\noindent Now, if $u_0$ is chosen so large that $\frac12\leq r_2$ in $[u_{0},+\infty)$,  

\begin{equation} \label{trick} \frac{1}{V_2}-\frac{1}{V_1}=\frac{V_1^2-V_2^2}{V_1V_2(V_1+V_2)}\leq \frac{2^{2N-2}w}{V_1V_2(V_1+V_2)}. 
\end{equation} 

\noindent Integrating \eqref{step5} and using \eqref{keller}, it follows that for $u\ge u_{0}$,

$$  w(u)\leq w(u_{0})+ {C}(N)\int_{u_0}^u \left(\int_{u'}^{+\infty} \frac{w}{F^\frac32}du''\right)f\;du'.$$

\noindent Integrating by parts

$$ w(u)\leq w(u_{0})+ {C}(N)\left(F(u)\int_u^{+\infty} \frac{w}{F^\frac32}du' + \int_{u_0}^u \frac{w}{F^\frac12}du' \right). $$

\noindent Thanks to the Keller-Osserman condition \eqref{KO}, if $u_{0}$ is chosen large enough,

$$\int_{u_0}^u \frac{w}{F^\frac12}du'\leq w(u)\int_{u_0}^{+\infty} \frac {1}{\sqrt F}\le\frac1{2C(N)}w(u).$$

\noindent We have then obtained

\begin{equation} \label{weq}  w(u)\leq 2w(u_{0})+ 2{C}(N)F(u)\int_u^{+\infty} \frac{w}{F^\frac32}du'. \end{equation} 

\noindent Introduce $G(u)=\int_u^{+\infty}{w\over F^\frac 3 2}du'$. Thanks to \eqref{keller} and \eqref{KO}, we have $G(+\infty)~=~0$. In addition, letting $c=2C(N)$, \eqref{weq} can be rewritten as  

$$ -\frac{dG}{du}\leq \frac{2w(u_{0})}{F^{\frac 3 2}}+ \frac{c}{F^{\frac12}}G.$$
That is,

$$ -\frac{d}{du}\left(G\exp\left(-c\int_u^{+\infty} \frac{1}{\sqrt F}du'\right)\right)\leq \frac{2w(u_{0})}{F^{\frac 3 2}}\exp\left(-c\int_u^{+\infty} \frac{1}{\sqrt F}du'\right)\leq \frac{2w(u_{0})}{F^{\frac 3 2}}. $$

\noindent Integrating between $u$ and $+\infty$, we then obtain, using once again \eqref{KO},
  
$$ G(u)\leq C \int_u^{+\infty}  \frac{1}{F^{\frac 3 2}}= o\left(\frac 1 F\right).$$

\noindent Going back to \eqref{weq}, we deduce that $w$ is bounded above. 
%

\noindent {\bf Step 6.} The difference $U_{2}(r)-U_{1}(r)$ converges to $0$ as $r\to1$.

Given $r$ close to $1$ and $i\in\{1,2\}$, let $u_{i}=U_{i}(r)$. Then,
\begin{equation*} 
\int_{u_1}^{+\infty} \frac 1 {V_1}du=1-r =\int_{u_2}^{+\infty} \frac 1 {V_2}du.
\end{equation*}
That is,
$$ \int_{u_1}^{u_2} \frac 1 {V_1}du= \int_{u_2}^{+\infty} \left(\frac 1 {V_2}-\frac 1 {V_1}\right)du.
$$
\noindent Using \eqref{trick}, \eqref{keller}, and the previous step, we deduce that

$$ \int_{u_1}^{u_2} \frac 1 {\sqrt F}du\le C\int_{u_2}^{+\infty} \frac{1}{F^{3/2}}du.$$
It follows that
$$ 0\le \frac{u_2-u_1}{\sqrt{F(u_2)}}\leq \frac{C}{\sqrt{F(u_2)}}\int_{u_2}^{+\infty} \frac 1 F du$$

\noindent and the claim follows promptly.
}

\noindent{\bf Step 7.} End of proof.

Let $w=U_{2}-U_{1}$.  Since $U_{2}\ge U_{1}$ and $f$ is nondecreasing, we see from the previous step that
\begin{equation*}
 \left\{
 \begin{aligned}
\Delta w &=f(U_{2})-f(U_{1})\ge 0 &\quad\text{in $B$,}\\
w &= 0 &\quad\text{on $\partial B$.}
\end{aligned}
\right.
 \end{equation*}
By the maximum principle, $w\le 0$ in $B$, as desired.

\section{Proof of Theorem \ref{thm2}} 
Take a solution $u$ to \eqref{main}. Let $a$ be the constant appearing in \eqref{positivity}, $M$ the constant beyond which $\sqrt F$ is convex, and fix $\tilde M>\max(0,a,M)$.  Fix $\eps>0$ so small that $u>\tilde M$ in $\Omega_{\eps}=\{x\in\Omega : \dist(x,\partial\Omega)<\eps\}$.  

\noindent{\bf Step 1.}
We begin by proving that there exists a sequence of functions $(u_{N})_{N\in\N}$ solving
\begin{equation}\label{un} 
 \left\{
 \begin{aligned} 
\Delta u_{N} &= f(u_{N}) &\quad\text{in $\Omega_{\eps}$,}\\
u_{N} &= N &\quad\text{on $\partial \Omega$,}\\
u_{N} &= u &\quad\text{on $\{x\in\Omega : \dist(x,\partial\Omega)=\eps\}$,}
\end{aligned}
\right. 
 \end{equation} 
such that
\begin{equation} \label{esti} 
0\le u_{N}\le u\qquad\text{in $\Omega_{\eps}$.}
\end{equation}    
We may always assume that $f(0)=0$.\footnote{If $f(0)\neq0$, work with  any nondecreasing $C^1$ function $\tilde f$ such that $\tilde f(0)=0$, $\tilde f=f$ on $[\tilde M,+\infty)$.}
In particular, $\sub u=0$ and $\super u=u$ are respectively a sub and supersolution of \eqref{un} and they are ordered. It follows that there exists a solution $u_{N}$ to \eqref{un} such that \eqref{esti} holds. 

A standard application of the maximum principle shows that $u_{N}$ is the unique solution to \eqref{un} and that $(u_{N})$ is a nondecreasing sequence. Thanks to \eqref{esti} and elliptic regularity, we may also assert that $(u_{N})$ converges in $C^2_{loc}(\super\Omega_{\eps}\setminus\partial\Omega)$ to a function $\tilde u$ solving
\begin{equation}\label{tilde u} 
 \left\{
 \begin{aligned} 
\Delta \tilde u &= f(\tilde u) &\quad\text{in $\Omega_{\eps}$,}\\
\tilde u &= +\infty &\quad\text{on $\partial \Omega$,}\\
\tilde u &= u &\quad\text{on $\{x\in\Omega : \dist(x,\partial\Omega)=\eps\}$,}
\end{aligned}
\right. 
 \end{equation} 
 
\noindent{\bf Step 2.} 
There holds
\begin{equation} \label{esti2} 
\vert \nabla u_{N}\vert^2 - 2F(u_{N}) \le M_{N} \qquad\text{in $\Omega_{\eps}$}, 
\end{equation} 
where 
\begin{equation} \label{mn} 
M_{N} = \sup_{\dist(x,\partial\Omega)=\eps}\left[\vert \nabla u_{N}\vert^2 - 2F(u_{N})\right].
\end{equation}  
The proof is a straightforward adaptation of an argument due to Bandle and Marcus (\cite{bandle-marcus-jam}), which uses the method of $P$-functions. We give the full argument here for convenience of the reader.
Let 
$$P_{N} = \vert \nabla u_{N}\vert^2 - 2F(u_{N}).$$ 
By a result of Payne and Stackgold (\cite{PS}, see also Chapter 5 in \cite{S}), there exists a bounded continuous vector field $A$, such that 
$$
\Delta P_{N} - \frac{A\cdot\nabla P_{N}}{\vert \nabla u_{N}\vert^2}\ge0
$$
at every point in $\Omega_{\eps}$ where $\nabla u_{N}\neq0$. Hence, $P_{N}$ attains its maximum over $\super\Omega_{\eps}$ either on $\partial\Omega$, on $\{x\in\Omega : \dist(x,\partial\Omega)=\eps\}$, or at a critical point of $u_{N}$. It only remains to prove that the first case cannot occur. We claim that $\partial P_{N}/\partial n\le0$ on $\partial\Omega$, where $n$ is the outward unit normal to $\partial\Omega$. The boundary-point lemma then implies that $P_{N}$ cannot attain its maximum on $\partial\Omega$. It remains to prove our claim. Observe that since $u_{N}$ is constant on $\partial\Omega$, $\vert \nabla u_{N}\vert=\partial u_{N}/\partial n$ on $\partial\Omega$. Hence,
$$
\frac{\partial P_{N}}{\partial n} = 2 \frac{\partial u_{N}}{\partial n}\frac{\partial^2 u_{N}}{\partial n^2} - 2f(N)\frac{\partial u_{N}}{\partial n},\qquad\text{on $\partial\Omega$.}
$$
Furthermore, letting $H$ denote the mean curvature of $\partial\Omega$,
$$
\Delta u_{N} = \frac{\partial^2 u_{N}}{\partial n^2} + (N-1)H\frac{\partial u_{N}}{\partial n}\qquad\text{on $\partial\Omega$.}
$$
Since $\partial u_{N}/\partial n>0$ and $H\ge0$, this implies that
$$
\frac{\partial^2 u_{N}}{\partial n^2} -f(N)\le0
$$
and consequently $\partial P_{N}/\partial n\le0$, as desired. We have just proved \eqref{esti2}. 

\noindent{\bf Step 3.} The function $\tilde u=\lim_{N\to+\infty}u_{N}$ coincides with $u$ in $\Omega_{\eps}$. 

The proof of this fact bears resemblances with a trick due to L. Nirenberg given in \cite{brezis-kamin}. By \eqref{esti}, we already have $\tilde u\le u$ in $\Omega_{\eps}$ and it remains to prove the reverse inequality. Thanks to \eqref{esti} and elliptic regularity, there exists a constant $M$ such that
$$
2M\ge M_{N},
$$ 
where $M_{N}$ is given by \eqref{mn}. Now let $\tilde F = F +M$ and define 
$$
v_{N} = \int_{u_{N}}^{+\infty}\frac{dt}{\sqrt{2\tilde F(t)}}.
$$ 
Then, \eqref{esti2} can be rewritten as
$$
\vert \nabla v_{N}\vert \le 1\qquad\text{in $\Omega_{\eps}$}
$$ 
from which it easily follows that
\begin{equation} \label{gradinf} 
\vert \nabla \tilde v\vert \le 1\qquad\text{in $\Omega_{\eps}$},
\end{equation} 
where we defined similarly
$$
\tilde v = \int_{\tilde u}^{+\infty}\frac{dt}{\sqrt{2\tilde F(t)}}.
$$
Let at last
$$
v = \int_{u}^{+\infty}\frac{dt}{\sqrt{2\tilde F(t)}}.
$$
It remains to prove that $u\le\tilde u$ i.e. $\tilde v\le v$ in $\Omega_{\eps}$. Using the equations satisfied by $u$ and $\tilde u$, we see that $w=v-\tilde v$ solves
\begin{align*}
-\Delta w &= \frac{f}{\sqrt{2\tilde F}}(u) \left(1-\vert\nabla v\vert^2\right) - \frac{f}{\sqrt{2\tilde F}}(\tilde u) \left(1-\vert\nabla \tilde v\vert^2\right)\\
&= \left[\frac{f}{\sqrt{2\tilde F}}(u)-\frac{f}{\sqrt{2\tilde F}}(\tilde u)\right]\left(1-\vert\nabla \tilde v\vert^2\right) + \frac{f}{\sqrt{2\tilde F}}(u)\left(\vert\nabla \tilde v\vert^2-\vert\nabla  v\vert^2\right)
\end{align*}
Since $\sqrt{2F}$ is convex, $\frac{f}{\sqrt{2\tilde F}}$ is nondecreasing. Using this and \eqref{gradinf}, we deduce that 
\begin{equation*}
 \left\{
 \begin{aligned} 
-\Delta w + b(x)\cdot\nabla w&\ge0,&\quad\text{in $\Omega_{\eps}$}
\\
w &= 0 &\quad\text{on $\partial \Omega_{\eps}$,}
\end{aligned}
\right. 
\end{equation*}  
where $b(x)=\frac{f}{\sqrt{2\tilde F}}(u)\nabla(v+\tilde v)$ is locally bounded in $\Omega$. We may now apply the maximum principle to conclude that $w\ge0$ in $\Omega$, as desired.

\noindent{\bf Step 4.} End of proof. The rest of the proof is similar to an argument due to Garcia-Melian \cite{garmel}. We take two arbitrary solutions $u,\super u$ of our equation \eqref{main}. We let $u_{N}, \super u_{N}$ be the corresponding solutions to the approximated problem \eqref{un}. In particular, $w_{N}= u_{N}-\super u_{N}$ solves
\begin{equation}\label{wn} 
 \left\{
 \begin{aligned} 
\Delta w_{N} &=f(u_{N})-f(\super u_{N}) &\quad\text{in $\Omega_{\eps}$,}\\
w_{N} &= 0 &\quad\text{on $\partial \Omega$,}\\
w_{N} &= u-\super u &\quad\text{on $\{x\in\Omega : \dist(x,\partial\Omega)=\eps\}$,}
\end{aligned}
\right. 
\end{equation}  
By the maximum principle,
$$
w_{N}\le \sup_{\dist(x,\partial\Omega)=\eps}(u-\super u)\qquad\text{in $\Omega_{\eps}$},
$$
with equality at some point $x_{N}$ such that $\dist(x_{N},\partial\Omega)=\eps$. Extracting a sequence if necessary, we deduce that $w=u-\super u$ satisfies
\begin{equation} \label{esti3} 
w\le \sup_{\dist(x,\partial\Omega)=\eps}(u-\super u)\qquad\text{in $\Omega_{\eps}$},
\end{equation} 
with equality at some point $z$ such that $\dist(z,\partial\Omega)=\eps$. Now, we also have
\begin{equation*}
\left\{
 \begin{aligned} 
\Delta w &=f(u)-f(\super u) &\quad\text{in $\Omega\setminus\Omega_{\eps}$,}\\
w &= u-\super u &\quad\text{on $\{x\in\Omega : \dist(x,\partial\Omega)=\eps\}$.}
\end{aligned}
\right. 
\end{equation*}  
By the maximum principle, we deduce that inequality \eqref{esti3} holds throughout $\Omega$, with equality at the point $z$. The strong maximum principle implies that $w$ is equal to a constant $c$. Since $u,\super u$ solve \eqref{main}, we deduce that $f(u)=f(u+c)$, which is possible only if $c=0$. 
\hfill\qed

\section{Proof of Proposition \ref{pro1}}
We thank Laurent V\'eron (\cite{veron}) for the following proof.
Given $p\in(1,N/(N-2))$, $k\in \N$ and $\lambda>0$, we begin by solving
\begin{equation}\label{deltak} 
 \left\{
 \begin{aligned} 
-\Delta u + u^p&=\lambda\delta_{0} &\quad\text{in $B$,}\\
u &= k &\quad\text{on $\partial B$,}
\end{aligned}
\right. 
 \end{equation} 
Since $0$ is a subsolution, while a large constant multiple of the fundamental solution is a supersolution, we deduce from the method of sub and supersolution (see e.g. \cite{mp} for the appropriate statement) that there exists a solution $u=u_{k}$ to \eqref{deltak}. By the maximum principle, $u_{k}$ is the unique solution to \eqref{deltak}, and the sequence $(u_{k})$ is nondecreasing. Thanks to the Keller-Osserman estimate (see e.g. \cite{keller}), the sequence $(u_{k})$ is uniformly bounded on compact subsets of the punctured ball $B\setminus\{0\}$.  It follows from elliptic regularity that $u_{k}$ converges to a solution  $u=u_{\lambda}$ of
\begin{equation*}
 \left\{
 \begin{aligned} 
-\Delta u + u^p&=\lambda\delta_{0} &\quad\text{in $B$,}\\
u &= +\infty &\quad\text{on $\partial B$,}
\end{aligned}
\right. 
 \end{equation*} 
By the results of \cite{veron2}, $u_{\lambda}$ behaves like a constant multiple of the fundamental solution near the origin. In particular, each $u_{\lambda}$ is a large solution in the punctured ball. 

There exists yet another large solution. Simply note that for an appropriate constant $c=c(N,p)>0$, the function $u_{1}(x)=c\vert x\vert^{-2/(p-1)}$ solves $\Delta u =u^p$ in $\R^N\setminus\{0\}$. Let also $u_{2}$ be the unique solution to 
\begin{equation*}
 \left\{
 \begin{aligned} 
\Delta u &=u^p &\quad\text{in $B$,}\\
u &= +\infty &\quad\text{on $\partial B$,}
\end{aligned}
\right. 
 \end{equation*} 
Then, $\sub u = \max(u_{1},u_{2})$ and $\super u = u_{1}+u_{2}$ form an ordered pair of sub and supersolution 
to the equation in the punctured ball. The method of sub and supersolutions implies the existence of a new large solution $u_{\infty}$ which behaves like $c\vert x\vert^{-2/(p-1)}$ near the origin, hence distinct from $u_{\lambda}$.

Finally, observe that for the nonlinearity $f(u)=u^p$, if $u$ is a large solution and $\epsilon>0$, then $(1+\epsilon)u$ is a supersolution. From this, the classification of singularities both at the origin (see \cite{veron2}) and on the boundary (see e.g. \cite{bandle-marcus-jam}), and the maximum principle, it easily follows that the set of positive large solutions in the punctured ball is exactly $\{u_{\lambda}\}_{\lambda\in(0,+\infty]}$.

\

{\noindent\bf Acknowledgements.} 
O. C. was  supported in part by NSF grants DMS-0807266 and DMS-1108794. Any opinions, findings, conclusions or recommendations expressed in this material are those of the authors and do not necessarily reflect the views of the National Science Foundation. L.D. was supported in part by ERC grant 277749 EPSILON .

\end{document}